\documentclass[handcarry]{aiaa-tc}
\usepackage{placeins}
\usepackage{indentfirst}
\usepackage{multirow}

\usepackage{amssymb}
\usepackage{amsmath,latexsym}
\usepackage{url}
\usepackage{enumerate}
\usepackage{graphicx}

\usepackage{booktabs}

\usepackage{mathtools}

\usepackage{subcaption}
\usepackage{cleveref}
\usepackage{graphicx}
\usepackage{amsmath}
\usepackage[version=4]{mhchem}
\usepackage{siunitx}
\usepackage{longtable,tabularx}
\usepackage{amsmath}

\usepackage{examplep}
\usepackage[driverfallback=dvipdfm,
bookmarks=true,    
bookmarksopen=true,
bookmarksopenlevel=\maxdimen,
bookmarksdepth=10
]{hyperref}
\hypersetup{
pdfcenterwindow=true, 
pdfpagelabels=true,
pagebackref=true,            
hypertexnames=true,
plainpages=false,
unicode=false,     
pdftoolbar=true,                
    pdfmenubar=true,          
    pdffitwindow=false,         
    pdfstartview={FitH},        
    pdftitle={On Iterative Convergence for Damped Viscous Scheme},
    pdfauthor={Hiroaki Nishikawa, Yoshitaka Nakashima, Norihiko Watanabe},    
    pdfsubject={On Iterative Convergence for Damped Viscous Scheme},       
    pdfcreator={Hiroaki Nishikawa},   
    pdfproducer={Hiroaki Nishikawa}, 
    baseurl={http://www.hiroakinishikawa.com},
    pdfkeywords={}, 
    pdfnewwindow=true,     
    colorlinks=true,       
    linkcolor=black,       
    citecolor=black,       
    filecolor=black,       
    urlcolor=black,        
  breaklinks=true,
  hyperfigures=true,
  backref=true,
breaklinks=false, 
pdfpagelabels,      
pagebackref,        
hypertexnames=true, 
plainpages=false,   
naturalnames        
}
\usepackage[all]{hypcap}
\numberwithin{equation}{section}
\numberwithin{subsubsection}{subsection}
\numberwithin{subsection}{section}
\renewcommand{\thesection}{\arabic{section}}

\usepackage[utf8]{inputenc}
\usepackage{dashbox}

\renewcommand{\thefootnote}{\arabic{footnote}}

\title{Improved Wall-Normal Derivative Formulae for Anisotropic Adaptive Simplex-Element Grids}

 \author{ {Hiroaki Nishikawa}\thanks{Research Fellow
({hiro@nianet.org}),
 100 Exploration Way, Hampton, VA 23666 USA, Associate Fellow AIAA}\\
  {\normalsize\itshape 
{National Institute of Aerospace},
 Hampton, VA 23666, USA}%
}

\def\o6{\frac{1}{6}} 

\begin{document}


 \maketitle

\begin{abstract}
In this paper, we explore methods for computing wall-normal derivatives used for calculating wall skin friction and heat transfer 
over a solid wall in unstructured simplex-element (triangular/tetrahedral) grids generated by anisotropic grid adaptation. 
Simplex-element grids are considered as efficient and suitable for automatic grid generation and adaptation, but present a challenge
to accurately predict wall-normal derivatives. For example, wall-normal derivatives computed by a simple finite-difference approximation, 
as typically done in practical fluid-dynamics simulation codes, are often contaminated with numerical noise. 
To address this issue, we propose an improved method based on a common step-length for the finite-difference approximation, which
is otherwise random due to grid irregularity and thus expected to smooth the wall-normal derivative distribution over a boundary. 
Also, we consider using least-squares gradients to compute the wall-normal derivatives and discuss their possible improvements. 
Numerical results show that the improved methods greatly reduce the noise in the wall-normal derivatives for irregular simplex-element grids. 
\end{abstract}
 
\section{Introduction}
\label{intro} 

Towards fully automated computational fluid dynamics (CFD) simulations, there is an increasing interest in developing robust and 
accurate discretizations and solvers for unstructured simplex-element grids with anisotropic grid adaptation [\citen{WhiteNishikawaBaurle_scitech2020,Kleb_etal_aiaa2019-2948,Nishikawa_scitech2020,nishikawa_centroid:JCP2020,uns_grid_adaptation_aiaa2018-1103,alauzet-loseille-decade-aniso-adapt-cfd}]. While fully-automated CFD simulations have recently been successfully demonstrated for practical three-dimensional problems [\citen{Kleb_etal_aiaa2019-2948}], improvements are still desired in the discretization method especially in terms of accuracy. For example, computations of derivatives such as the viscous/heat fluxes and the vorticity are known to suffer from numerical noise on highly irregular and skewed simplex-element grids [\citen{liu_nishikawa_aiaa2016-3969}]. Impact of numerical noise may be reduced in integrated quantities (e.g., drag coefficient), but when accuracy is required for a pointwise quantity, e.g., heating rate at a stagnation point, numerical noise can be a serious problem, causing a large error and a great degree of uncertainty in the prediction. In order to develop a reliable automated CFD solver, it is necessary to be able to accurately predict pointwise derivatives, especially at wall boundaries, on irregular anisotropic simplex-element grids. 

There have been efforts to improve derivative accuracy on unstructured grids: derivative variables introduced in target equations [\citen{nishikawa_hyperbolic_poisson:jcp2020,liu_nishikawa_aiaa2016-3969}], improved least-squares methods [\citen{NovelGradStencil:CF2018,SozerChristophCetin:AIAA2014,ShimaKitamuraHaga_AIAAJ2013}], and implicit gradient 
methods [\citen{WangRenPanLi:JCP2017,nishikawa_igg:JCP2018}]. All these methods are promising for improving derivative accuracy,
but they are not specifically designed to accurately estimate derivatives at a boundary and thus further improvements may be achieved by 
developing a more specific technique. In this work, we focus on wall-normal derivatives and explore various methods for 
accurately estimating the wall-normal derivatives relevant to the viscous stresses
and heat fluxes at a wall, as a post-processing step for a given numerical solution on adaptive simplex-element grids. A popular technique used in practical 
CFD codes for computing the wall-normal derivatives is 
a one-sided finite-difference formula applied in the wall-normal direction with the numerical solution in a nearby node/cell and a boundary condition value.
For example, it provides an approximation to $\partial u/ \partial n$ that dominates the viscous stress at a wall, where $u$ is a wall-parallel velocity components and $n$ is the 
coordinate normal to the wall. It can provide an accurate approximation on structured grids, but will introduce numerical noise on irregular simplex-element grids [\citen{WhiteNishikawaBaurle_scitech2020}]. 
We conjecture that the numerical noise in the derivative is due to the irregularity of the spacing used in the finite-difference formula. That is, the distance between the center of a cell
adjacent to a boundary and a wall-boundary face centroid varies along a boundary in adaptive simplex-element grids. Then, we propose to keep the same distance for the finite-difference 
formula applied at all wall-boundary faces. Numerical experiments show that the modified finite-difference method greatly improves the skin friction distribution over a flat plate
in a laminar flow, where the distance is determined as $\eta=0.5$ at all boundary faces, where $\eta$ is a normalized coordinate in the Blasius boundary-layer solution. 
In order to extend this technique to more general flows over complex geometries, it is necessary to develop a method for estimating a universal distance for the finite-difference formula.

Another post-processing technique can be devised based on gradients computed by a least-squares method at nodes or cells. Following a
recent study [\citen{WhiteNishikawaBaurle_scitech2020}], where heating rate prediction over a flat plate in a hypersonic turbulent flow
has been shown to be greatly improved by the face-averaged nodal-gradient (F-ANG) approach applied to the finite-volume discretization [\citen{WhiteNishikawaBaurle_scitech2020,NishikawaWhite_FANG:jcp2020}], we employ the F-ANG approach in this study, but only consider the post-processing technique for estimating the wall-normal derivatives for a given numerical solution. Discussions on the impact of the gradient approach (used in the discretization) on wall heating prediciton are given in a companion paper [\citen{NishikawaWhite_scitech2021-xxxx}]. 
 In the F-ANG method, gradients are computed
and stored at nodes while numerical solutions are stored at cells. Therefore, the wall-normal derivatives are already available at wall nodes.
However, these derivatives can be noisy for irregular anisotropic simplex-element grids. A simple method for reducing the noise would be to average the nodal 
gradients, and there are two possible techniques: face and cell averages. We compare these averaged gradients for estimating the wall-normal derivatives.
Numerical results indicate that the cell-averaged nodal-gradients can greatly reduce the noise in the wall-normal derivatives, but not as much as the improved 
finite-difference formula. 


\section{Finite-Volume Discretization}
\label{discretization} 

\begin{figure}[h!]
\begin{center}
\begin{minipage}[b]{0.67\textwidth}
\begin{center}
        \includegraphics[width=0.9\textwidth,trim=0 0 0 0,clip]{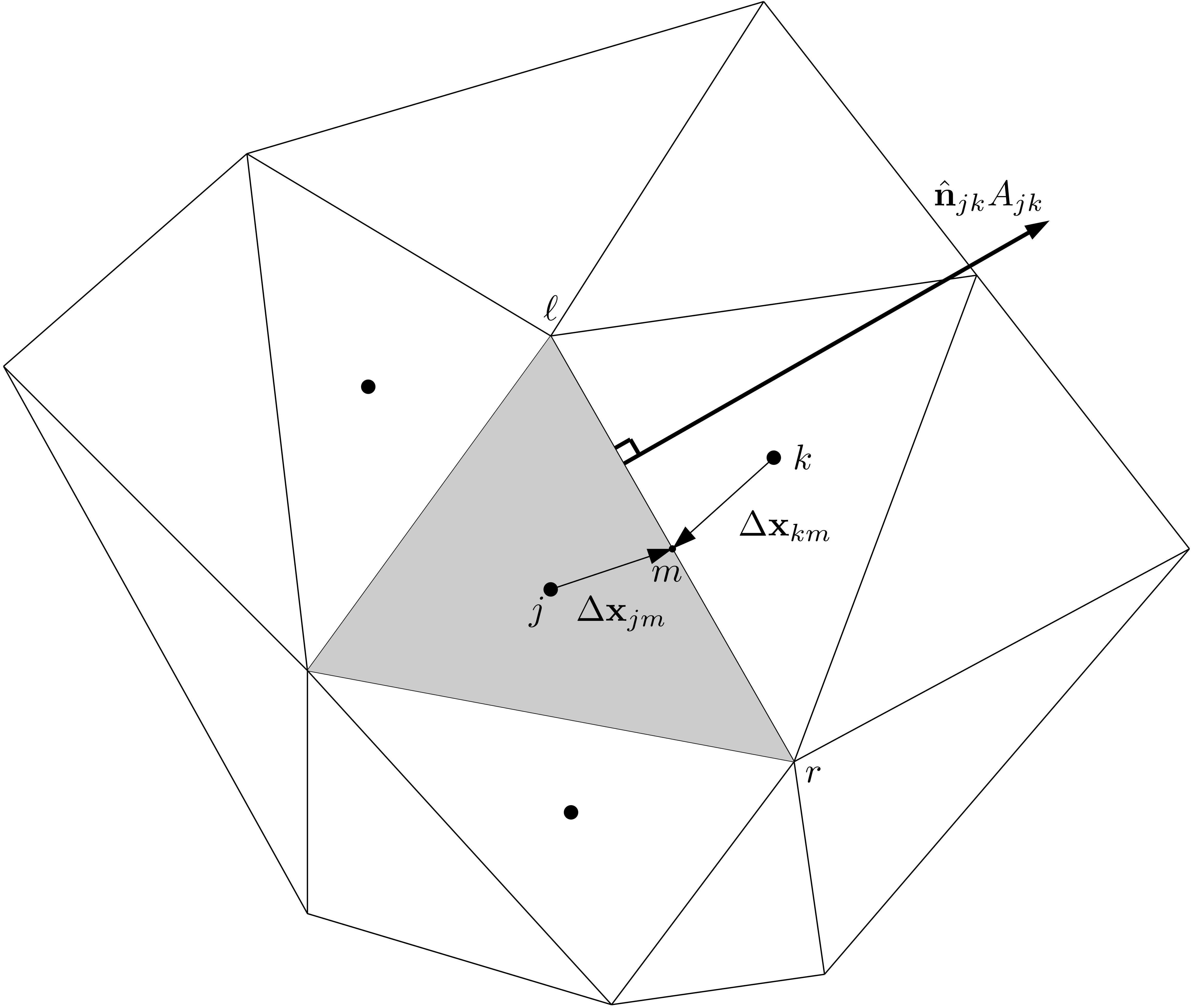}
\caption{Stencil for cell-centered finite-volume discretization.}
\label{fig:ccfv_stencil}
\end{center}
\end{minipage}
\end{center}
\end{figure}

We consider a cell-centered finite-volume discretization for the Navier-Stokes equations integrated 
over a computational cell $j$ (see Figure \ref{fig:ccfv_stencil}) with the midpoint rule: 
\begin{eqnarray}
{\bf Res}_j = \sum_{k \in \{ k_j \}}  {\bf \Phi}_{jk} A_{jk} ,
\label{res_used}
\end{eqnarray} 
where $ \{ k_j \}$ is a set of neighbors of the cell $j$, $ A_{jk}$ is the length of the face across $j$ and $k$, 
and ${\bf \Phi}_{jk} $ is a numerical flux. In this work, the Roe [\citen{Roe_JCP_1981}] flux is used for the inviscid terms, and
the alpha-damping flux [\citen{nishikawa:AIAA2010}] for the viscous terms. 
The solution values are stored in the primitive variables ${\bf w} = (\rho, {\bf v}, T)$, at cells as point values at the centroid [\citen{nishikawa_centroid:JCP2020}],
where $\rho$ is the density, ${\bf v}$ is the velocity vector, and $T$ is the temperature. 
Second-order accuracy is achieved by the linearly-exact flux quadrature (i.e., the midpoint rule) and 
the numerical flux evaluated with linearly reconstructed solutions. In this work, we focus on the F-ANG approach [\citen{WhiteNishikawaBaurle_scitech2020,NishikawaWhite_FANG:jcp2020}], which has
been shown to yield significantly more accurate wall-derivatives than a conventional finite-volume scheme, but still involves 
mild oscillations in the derivatives on irregular simplex-element grids. In the F-ANG approach, we compute the solution gradients
at nodes using solution values stored at cells by a linear least-squares method. Then, the nodal gradients are used in the linear reconstruction, averaged
over the nodes of a face: 
\begin{eqnarray}
  {\bf w}_L =   {\bf w}_j + \nabla {\bf w}_f \cdot \Delta {\bf x}_{jm}, \quad
  {\bf w}_R =   {\bf w}_k + \nabla {\bf w}_f \cdot \Delta {\bf x}_{km}, 
\label{states_LR}
\end{eqnarray} 
where $\nabla {\bf w}_f$ denotes the face-averaged gradient:
\begin{eqnarray}
\nabla {\bf w}_f  = \frac{1}{2} \left( 
\nabla {\bf w}_{jk}^\ell + \nabla {\bf w}_{jk}^r
\right),
\end{eqnarray}
where the superscripts $\ell$ and $r$ denote the left and right nodes of the face $[j,k]$ (see Figure \ref{fig:ccfv_stencil}). 
 This approach has various advantages especially for
simplex-element grids: e.g., no interpartition is required for gradients, the residual stencil is greatly reduced. See [\citen{WhiteNishikawaBaurle_scitech2020,NishikawaWhite_FANG:jcp2020}]
for further details.

\section{Methods for Computing Wall-Normal Derivatives}
\label{methods_wall_normal} 
 
In this section, we describe methods for computing a wall-normal derivative, taking the velocity derivative $\partial u/\partial n$
as an example, where $u$ is the wall-parallel component of the velocity and $n$ denotes the coordinate normal to the wall. 
We assume that the value at a wall $u_b$ is known, e.g., $u_b=0$ for a no-slip condition. The wall-normal derivative will be computed at a boundary node
or at a boundary face, depending on the method.

\begin{figure}[h!]
\begin{center}
\begin{minipage}[b]{0.9\textwidth}
\begin{center}
        \includegraphics[width=0.7\textwidth,trim=0 0 0 0,clip]{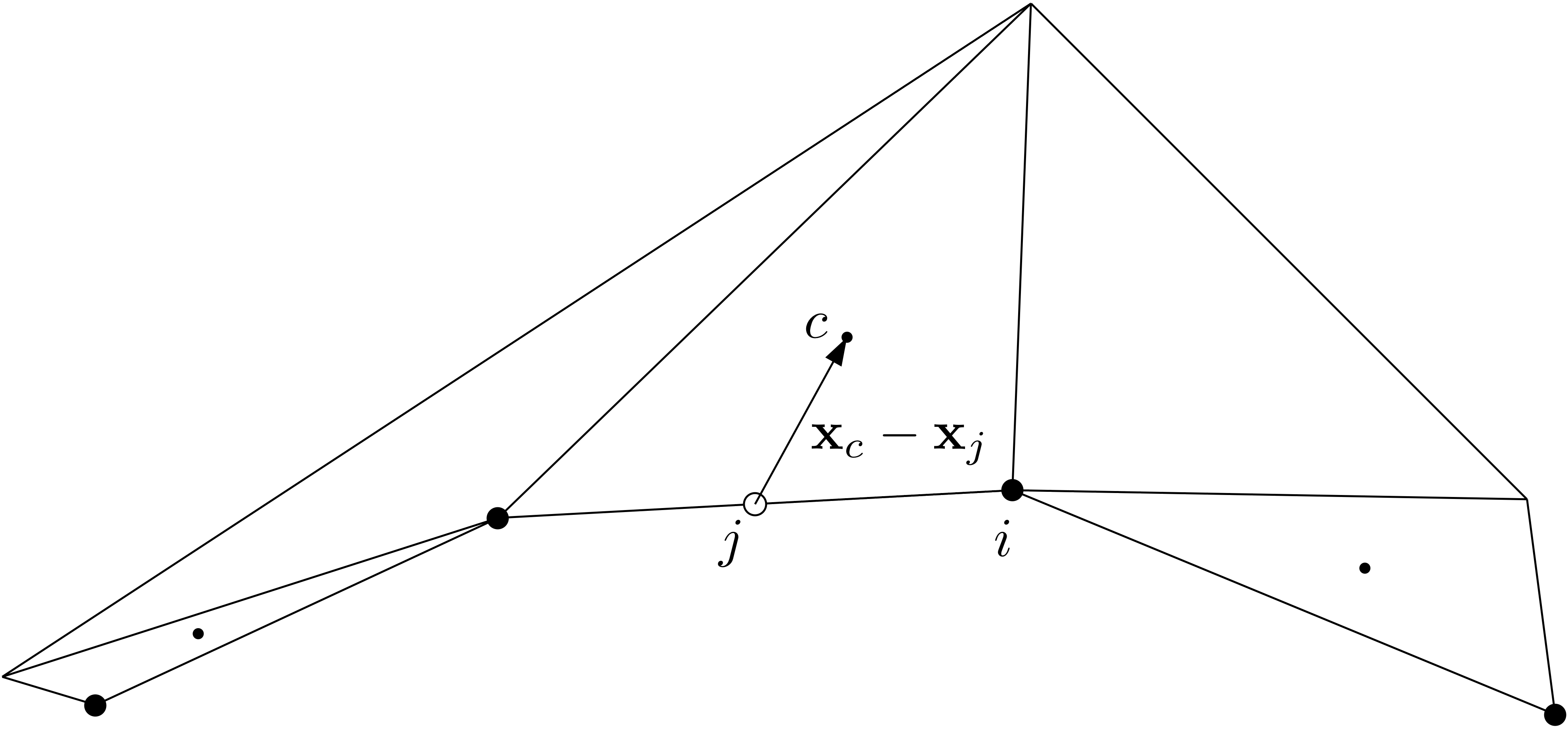}
\caption{Geometry of an irregular triangular grid over a boundary. Black filled circles are boundary nodes; smaller circles are
the centroids of triangles adjacent to the boundary. }
\label{fig:ccfv_stencil_boundary}
\end{center}
\end{minipage}
\end{center}
\end{figure}

\subsection{Nodal gradients}
\label{wall_derivative_ng} 

The wall-normal derivative is computed at a boundary node $i$ by projecting the nodal gradient $ \nabla u_i$ 
in the direction of wall normal (see Figure \ref{fig:ccfv_stencil_boundary}): 
\begin{eqnarray}
\left. \frac{\partial u}{\partial n} \right|_i  =  \nabla u_i \cdot \hat{\bf n}_i,
\end{eqnarray}
where $\hat{\bf n}_i$ denotes a unit vector at the node $i$ normal to the wall pointing towards the interior domain.
This is perhaps the simplest method since the nodal gradients are already available in the F-ANG method.

\subsection{Face-averaged nodal gradients}
\label{wall_derivative_fang} 

The wall-normal derivative is computed at a face $j$ with a face-averaged gradient: 
\begin{eqnarray}
\left. \frac{\partial u}{\partial n} \right|_j  =  \overline{\nabla} u_j^{face} \cdot \hat{\bf n}_j,
\end{eqnarray}
where $\overline{\nabla} u_j^{face}$ is the average of the nodal gradients over the boundary face $j$. 
In Figure \ref{fig:ccfv_stencil_boundary}, $\overline{\nabla} u_j^{face}$ is obtained by averaging the nodal gradients at the left
and right nodes of the boundary face $j$.

\subsection{Cell-averaged nodal gradients}
\label{wall_derivative_cang} 

The wall-normal derivative is computed at a face $j$ with a cell-averaged gradient: 
\begin{eqnarray}
\left. \frac{\partial u}{\partial n} \right|_j  =  \overline{\nabla} u_j^{cell} \cdot \hat{\bf n}_j,
\end{eqnarray}
where  $\overline{\nabla} u_j^{cell}$ is the average of the nodal gradients over a cell adjacent to the face $j$. 
In Figure \ref{fig:ccfv_stencil_boundary}, $\overline{\nabla} u_j^{cell}$ is obtained by averaging the nodal gradients at 
the three nodes of the triangular cell $c$.

  \begin{figure}[htbp!]
    \centering
          \begin{subfigure}[t]{0.32\textwidth}
        \includegraphics[width=\textwidth,trim=0 0 0 0 ,clip]{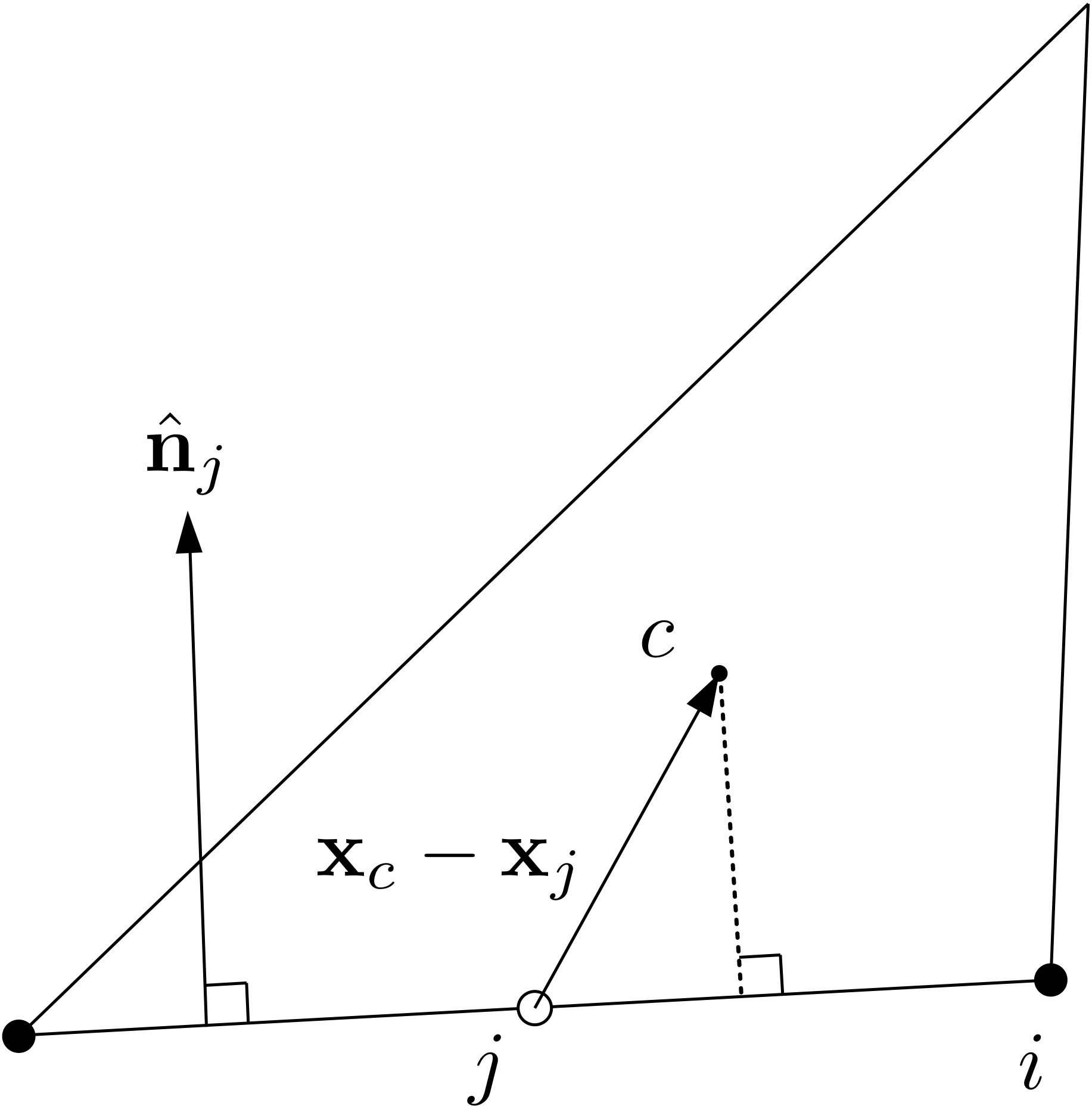}
          \caption{Finite-difference I.}
          \label{fig:geom_fd_megthod_I}
      \end{subfigure}
          \begin{subfigure}[t]{0.32\textwidth}
        \includegraphics[width=\textwidth,trim=0 0 0 0 ,clip]{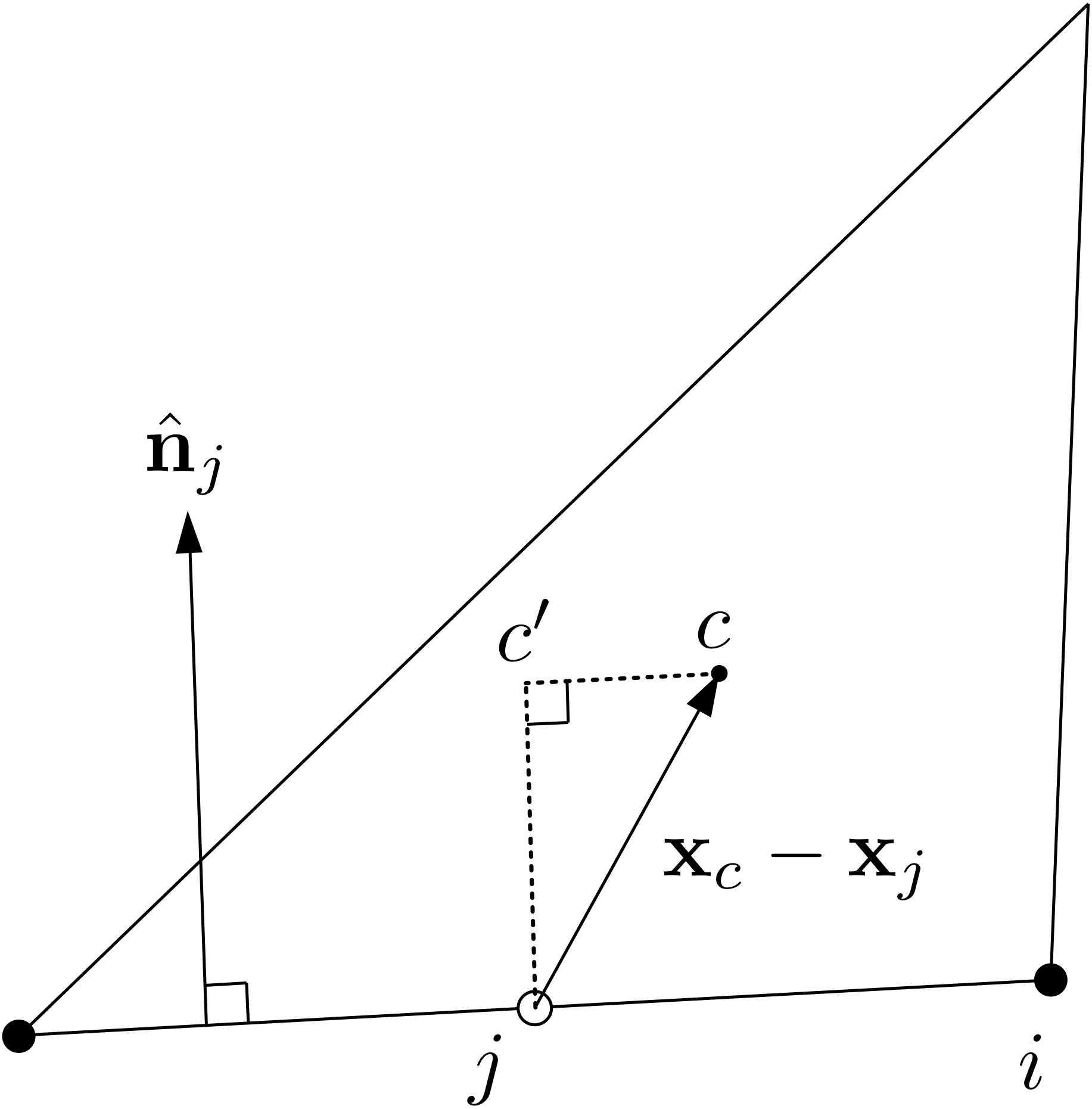}
          \caption{Finite-difference II.}
          \label{fig:geom_fd_megthod_II}
      \end{subfigure}
      \begin{subfigure}[t]{0.32\textwidth}
        \includegraphics[width=\textwidth,trim=0 0 0 0 ,clip]{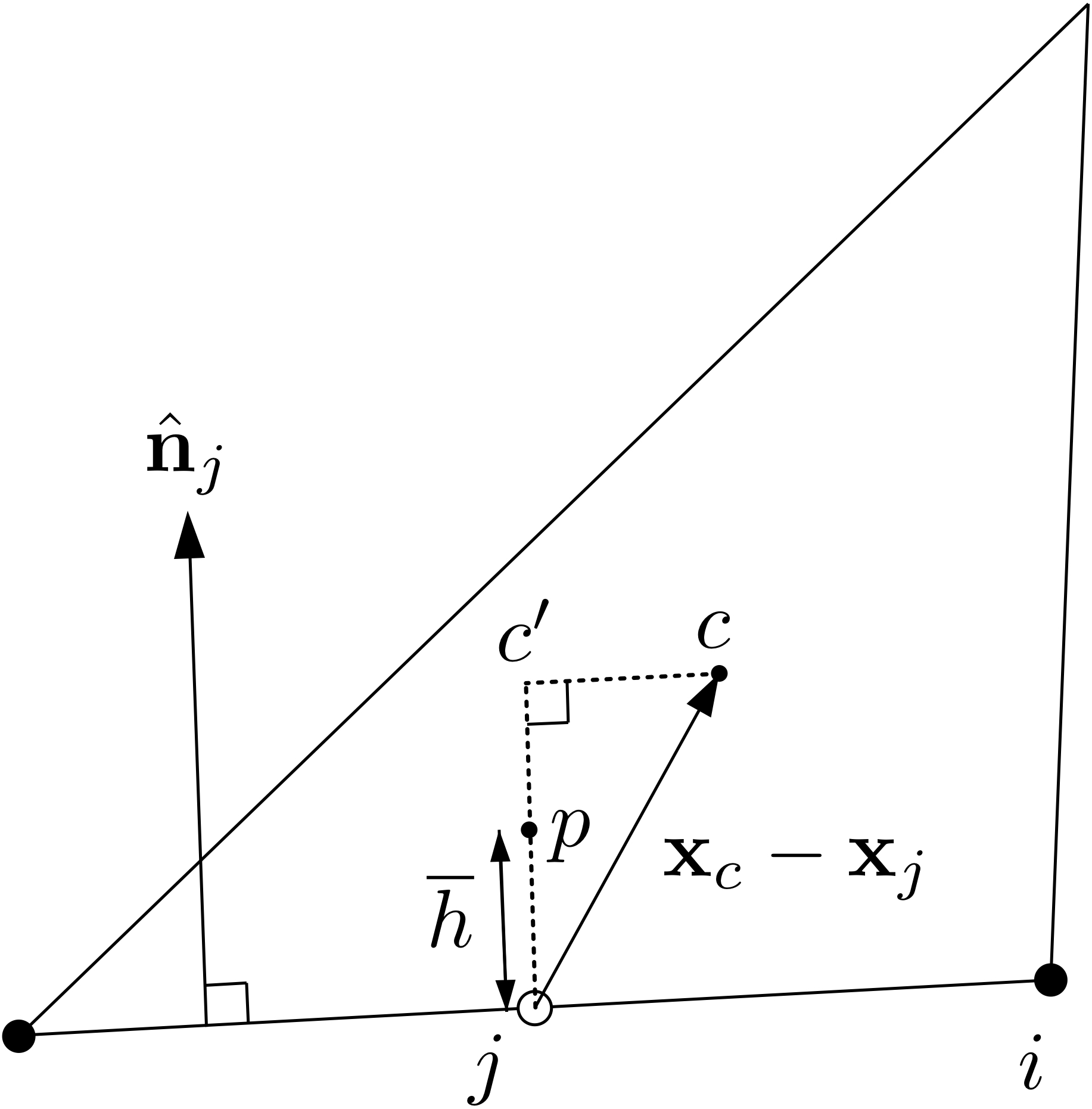}
          \caption{Finite-difference III.}
          \label{fig:geom_fd_megthod_III}
      \end{subfigure}
\caption{Geometries for the finite-difference wall-normal derivative formulae.}
\label{fig:geom_fd_megthod} 
\end{figure}

\subsection{Finite-difference I}
\label{wall_derivative_fd_1} 
 
The wall-normal derivative can be computed at a face $j$ by the following finite-difference formula:
\begin{eqnarray}
\left. \frac{\partial u}{\partial n} \right|_j  = 
\frac{ u_c  - u_b }{ ( {\bf x}_c - {\bf x}_ j ) \cdot \hat{\bf n}_j  },
\end{eqnarray}
where $u_c$ is a solution value stored at the cell $c$ adjacent to the boundary face $j$, 
$u_b$ is a known boundary value (e.g., from a boundary condition), 
${\bf x}_c$ is the centroid position of the cell $c$, $ {\bf x}_ j$ is the centroid
of the boundary face $j$. See Figure \ref{fig:geom_fd_megthod_I}.
This formula estimates the derivative 
at a point where the cell-center ${\bf x}_c$ is projected on the boundary, not at $ {\bf x}_ j$.

\subsection{Finite-difference II}
\label{wall_derivative_fd_2} 

The wall-normal derivative can be estimated at the face centroid by 
\begin{eqnarray}
\left. \frac{\partial u}{\partial n} \right|_j  = 
\frac{ u_{c'}  - u_b }{ ( {\bf x}_c - {\bf x}_ j ) \cdot \hat{\bf n}_j  },
\end{eqnarray}
where $u_{c'}$ is a solution extrapolated from the cell center ${\bf x}_c$ to the point
right above the face centroid (see Figure \ref{fig:geom_fd_megthod_II}): 
\begin{eqnarray}
u_{c'} = u_c + \nabla u_c \cdot ( {\bf x}_{c'}   -  {\bf x}_c   ), \quad
{\bf x}_{c'} = {\bf x}_{j} + [ ( {\bf x}_c - {\bf x}_ j ) \cdot \hat{\bf n}_j  ]  \hat{\bf n}_j.
\end{eqnarray}
This approximation is exact for linear functions at the face centroid.

\subsection{Finite-difference III}
\label{wall_derivative_fd_3} 

In the previous methods, the finite-difference formula is applied with varying step-length as the distance between
the cell center and the face center in the normal direction changes randomly on irregular simplex-element grids.  Hence, the methods are strongly dependent on 
the grid geometry, which is not suitable for adaptive simplex-element grids. 
To minimize the dependence on the grid geometry, we generalize the previous formula as
\begin{eqnarray}
\left. \frac{\partial u}{\partial n} \right|_j  = 
\frac{ u_p  - u_b }{ \overline{h} },
\end{eqnarray}
where $u_p$ is a solution extrapolated from the cell center ${\bf x}_c$ to the point 
at a certain distance $\overline{h}$ from the face center in the wall normal direction (see Figure \ref{fig:geom_fd_megthod_III}):
\begin{eqnarray}
u_p = u_c + \nabla u_c \cdot ( {\bf x}_p   -  {\bf x}_c   ), \quad
{\bf x}_p = {\bf x}_{j} +  \overline{h}   \, \hat{\bf n}_j,
\end{eqnarray}
where $\overline{h}$ is a global constant or varies along a wall but constant in a transformed normal coordinate.
This approximation is also exact for linear functions at the face centroid.

This formula does not depend strongly on the cell centroids of the adjacent cells. Whether in the physical or transformed coordinate, the above formula keeps the same step-length of the finite-difference formula over a wall and thus is expected to produce a smoother derivative (and its error) distribution.

In this work, we investigate the effects of the use of a common step-length $\overline{h}$ and explore estimation techniques of $\overline{h}$ 
for various viscous-flow problems on adaptive grids. Here, we provide preliminary results for a laminar flat-plate computation 
with the height $\overline{h}$ determined as described below.

\subsection{Finite-difference with face-area-weighted centroids}
\label{wall_derivative_with_remarkable_cetroids}
 
 Finite-difference III requires a common length to be defined. For a laminar flow over a flat plate, a normalized wall-normal distance, as we will describe later, has been found to serve well for subsonic flows. However, such an estimate is not applicable to high-speed flows and practical turbulent flows over complex geometries. As a practical technique, we consider the face-area-weighted centroid [\citen{nishikawa_centroid:JCP2020}], which could reduce the centroid-height variation along a boundary. The resulting method is similar to the finite-difference II but ${\bf x}_{c'}$ is replaced by the face-area-weighted centroid: \begin{eqnarray}{\bf x}_{c'} =  \frac{\displaystyle \sum_{k \in \{ k_j \}} \hat{A}_{jk}^p  {\bf x}_{jk}  }{\displaystyle  \sum_{k \in \{ k_j \}}  \hat{A}_{jk}^p}, \quad \hat{A}_{jk} = \frac{ {A}_{jk} }{ \displaystyle \max_{ k \in \{ k_j \}} {A}_{jk} }, \label{face_area_formula}\end{eqnarray}where${A}_{jk}$ is the area (length in two dimensions) of the face across $j$ and $k$,  ${\bf x}_{jk}$ is the face centroid, and $p=2$.

\section{Numerical Results}
\label{results}

\subsection{Two-dimensional laminar flat plate at $M_\infty = 0.15$, $Re_\infty = 10^6$}
\label{laminar_fp_mach0p15} 

For a laminar flat plate, we consider the normalized wall-normal coordinate $\eta$: 
\begin{eqnarray}
\eta = \frac{y}{x} \sqrt{Re_x},  \quad Re_x = \frac{\rho_\infty U_\infty x}{ \mu_\infty },
\end{eqnarray}
where $x$ is the distance from the leading edge, $\rho_\infty$, $U_\infty$, $\mu_\infty$ are free stream
density, flow speed, and viscosity. For the finite-difference III, we fix the height to be $\eta=0.5$, which is translated into
the $y$-coordinate as
\begin{eqnarray}
 \overline{h} = \eta \sqrt{ \frac{x_j}{Re_\infty} },
\end{eqnarray}
where $\eta=0.5$, $x_j$ is the $x$-coordinate of the face centroid, and $Re_\infty$ is the free stream Reynolds number 
per unit length. Note that we compute the wall-normal velocity derivative at boundary face centroids and thus $x_j$ will not
be zero. As can be seen from the formula, the height $ \overline{h}$ changes in the physical coordinate parabolically but 
remains constant in the normalized wall-normal coordinate $\eta$. 
 In the results, this method is denoted by FD-III($\eta=0.5$). For a more general flow problem, this approach may not be suitable.
Therefore, we also consider using a fixed length over the wall in the physical space: 
\begin{eqnarray}
 \overline{h} = \eta \sqrt{{Re_\infty} },
\end{eqnarray}
which corresponds to using the estimated distance with $\eta=0.5$ at $x_j=1.0$ everywhere. In this case, $ \overline{h} $ is a global constant.

\begin{figure}[htbp!]
    \centering
      \begin{subfigure}[t]{0.32\textwidth}
        \includegraphics[width=\textwidth,trim=4 4 4 4,clip]{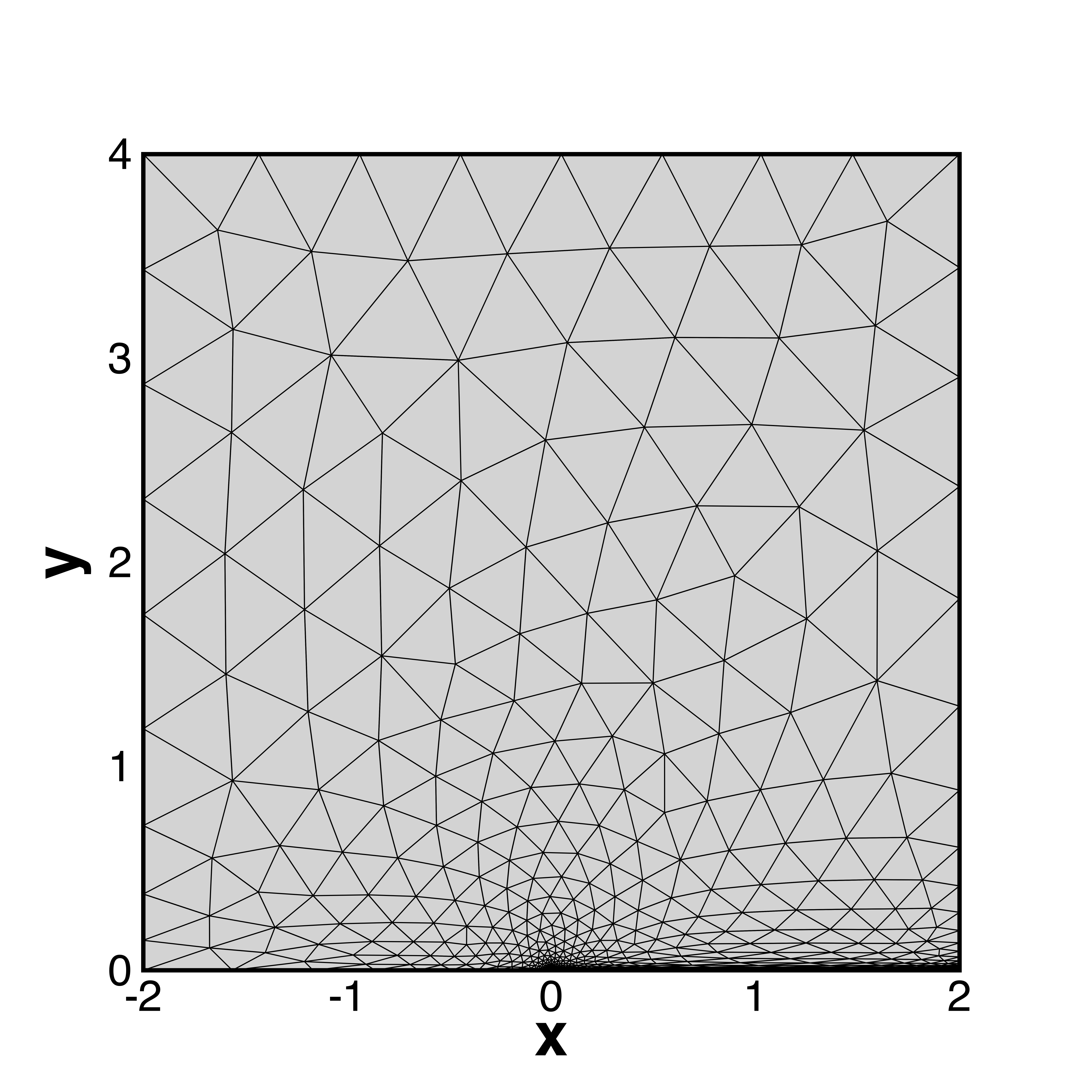}
          \caption{Adaptive grid (3457 nodes).}
          \label{fig:twod_laminar_fp_re1000000_grid}
      \end{subfigure}
          \begin{subfigure}[t]{0.32\textwidth}
        \includegraphics[width=\textwidth,trim=4 4 4 4,clip]{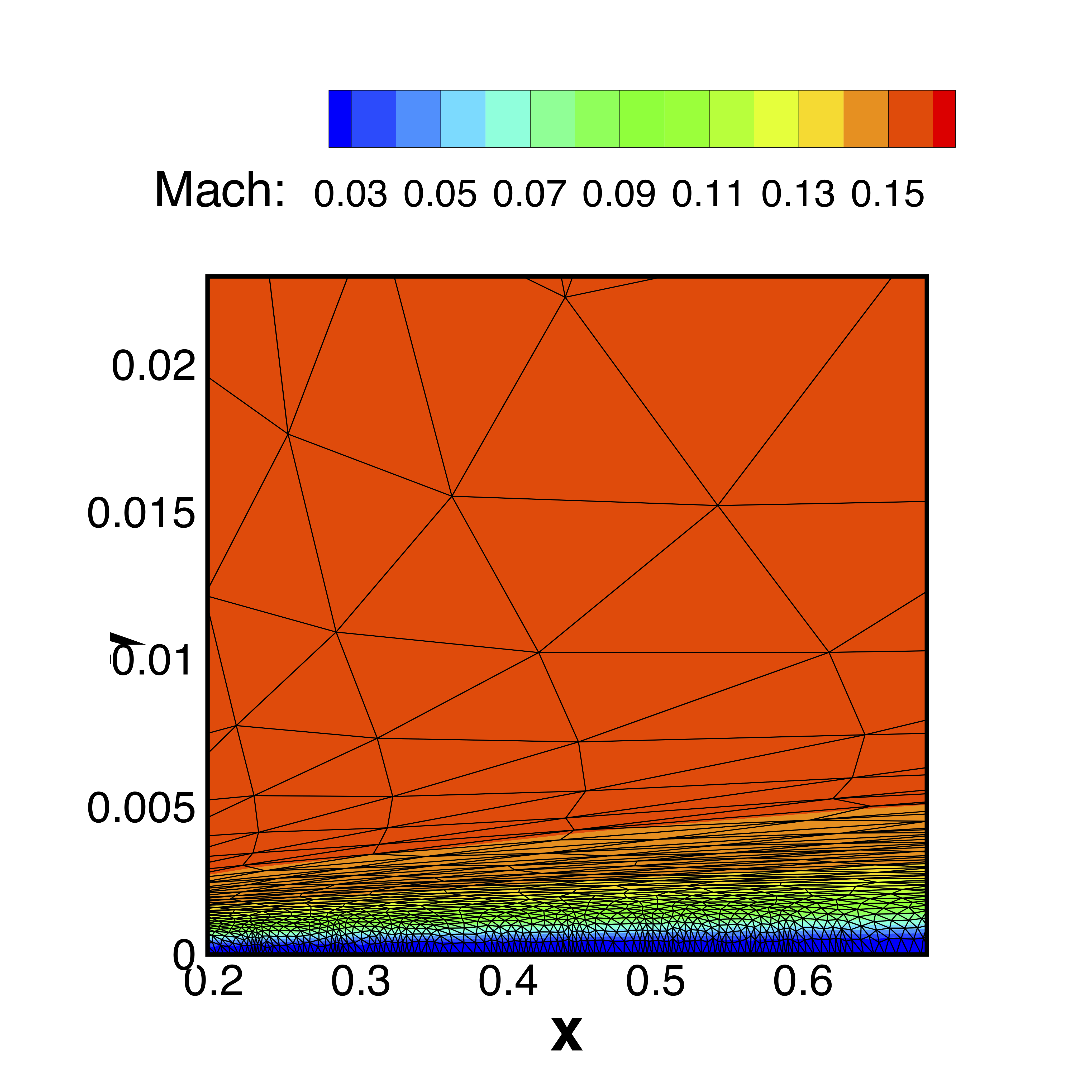}
          \caption{Zoomed-in view of the grid (with a scaled $y$-coordinate) and Mach contours.}
          \label{fig:twod_laminar_fp_re1000000_grid_zoom}
      \end{subfigure}
      \begin{subfigure}[t]{0.32\textwidth}
        \includegraphics[width=\textwidth,trim=4 4 4 4,clip]{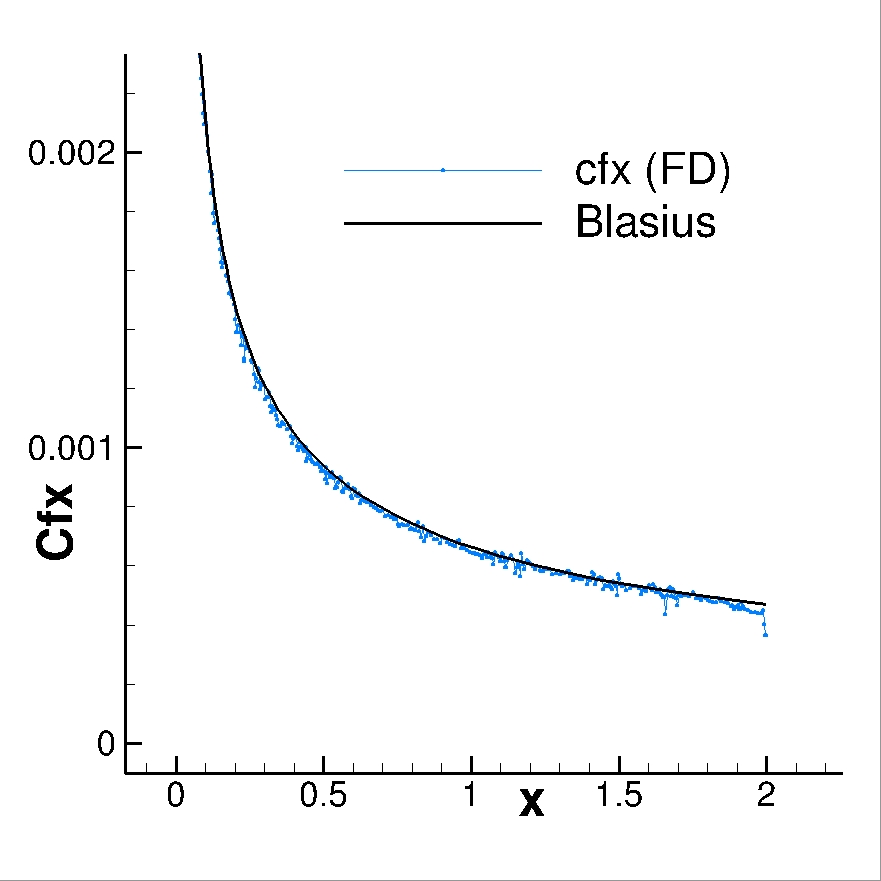}
          \caption{FD-I.}
          \label{fig:twod_laminar_fp_re1000000_grid_FD_I}
      \end{subfigure}
      \\
      \begin{subfigure}[t]{0.32\textwidth}
        \includegraphics[width=\textwidth,trim=4 4 4 4,clip]{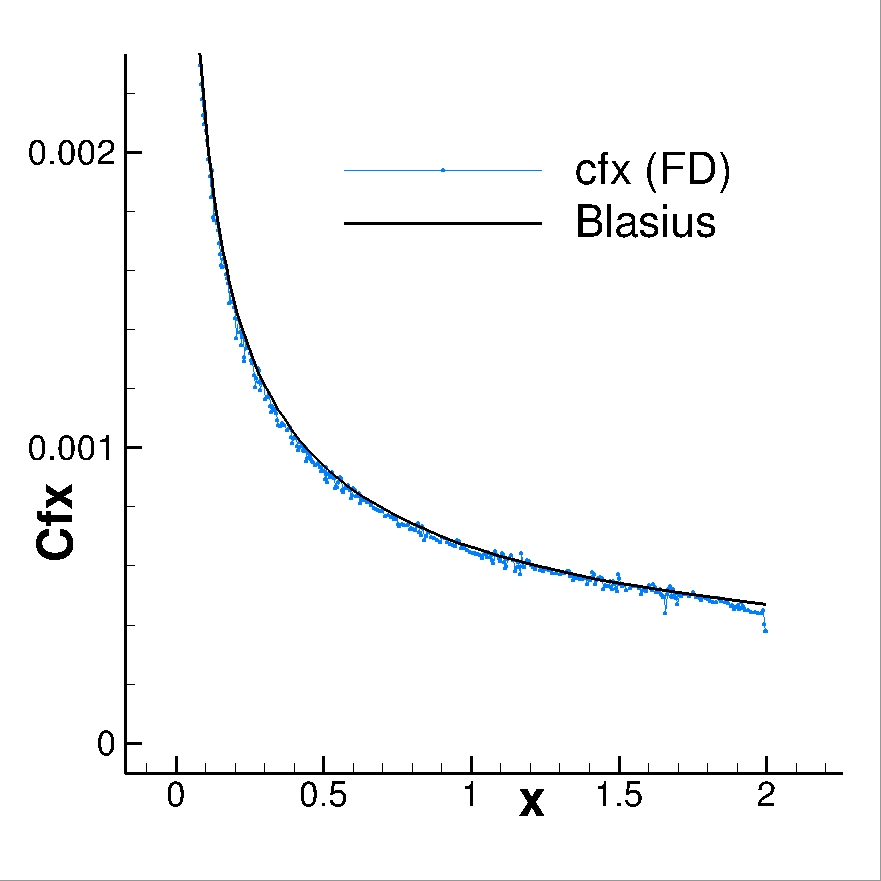}
          \caption{FD-II.}
          \label{fig:twod_laminar_fp_re1000000_grid_FD_II}
      \end{subfigure}
      \begin{subfigure}[t]{0.32\textwidth}
        \includegraphics[width=\textwidth,trim=4 4 4 4,clip]{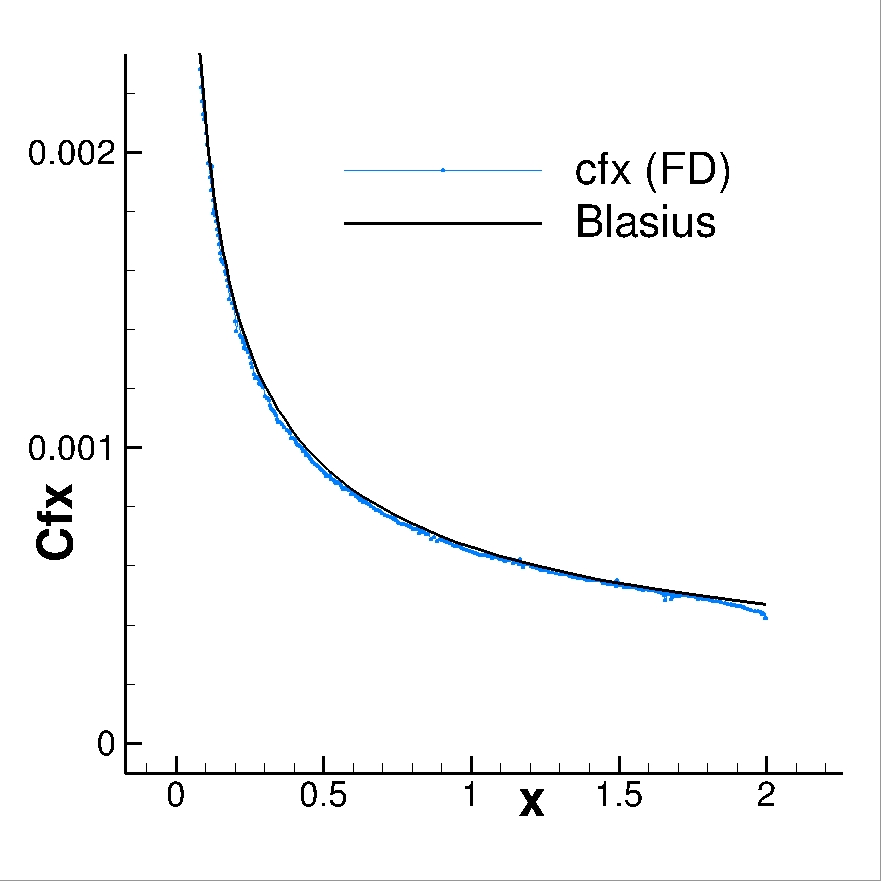}
          \caption{FD-III ($\eta=0.5$).}
          \label{fig:twod_laminar_fp_re1000000_grid_FD_III_eta}
      \end{subfigure}
      \begin{subfigure}[t]{0.32\textwidth}
        \includegraphics[width=\textwidth,trim=4 4 4 4,clip]{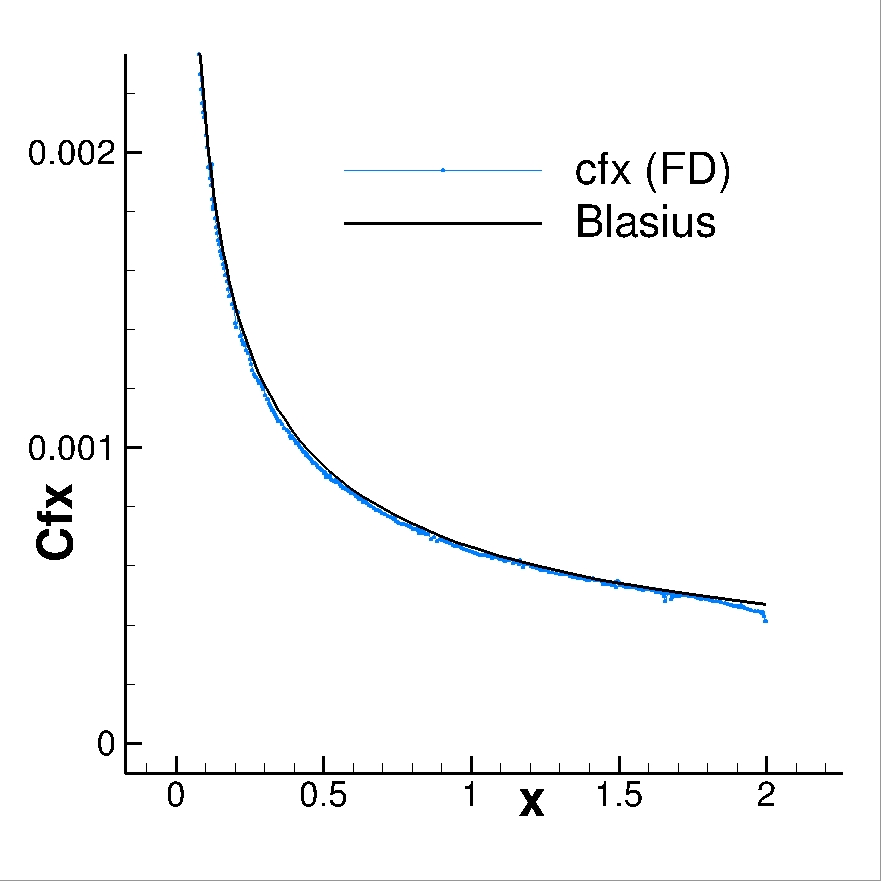}
          \caption{FD-III ($\eta=0.5$, $x=1$).}
          \label{fig:twod_laminar_fp_re1000000_grid_FD_III_eta_x}
      \end{subfigure}
      \\
      \begin{subfigure}[t]{0.32\textwidth}
        \includegraphics[width=\textwidth,trim=4 4 4 4,clip]{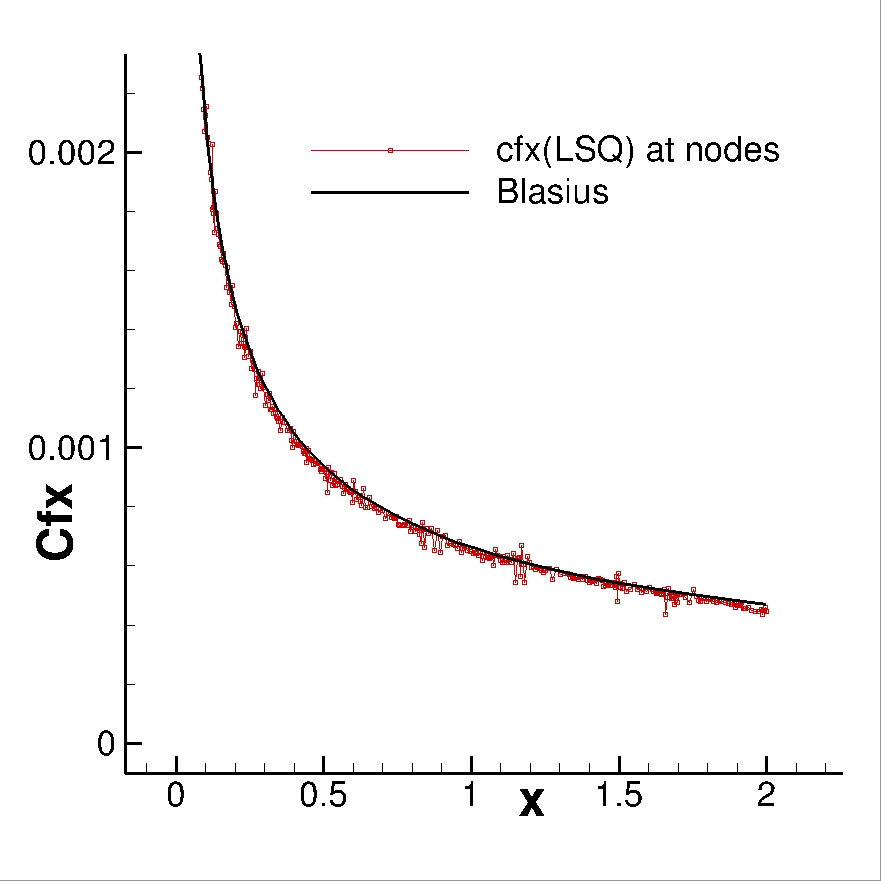}
          \caption{Nodal gradients.}
          \label{fig:twod_laminar_fp_re1000000_grid_LSQ}
      \end{subfigure}
      \begin{subfigure}[t]{0.32\textwidth}
        \includegraphics[width=\textwidth,trim=4 4 4 4,clip]{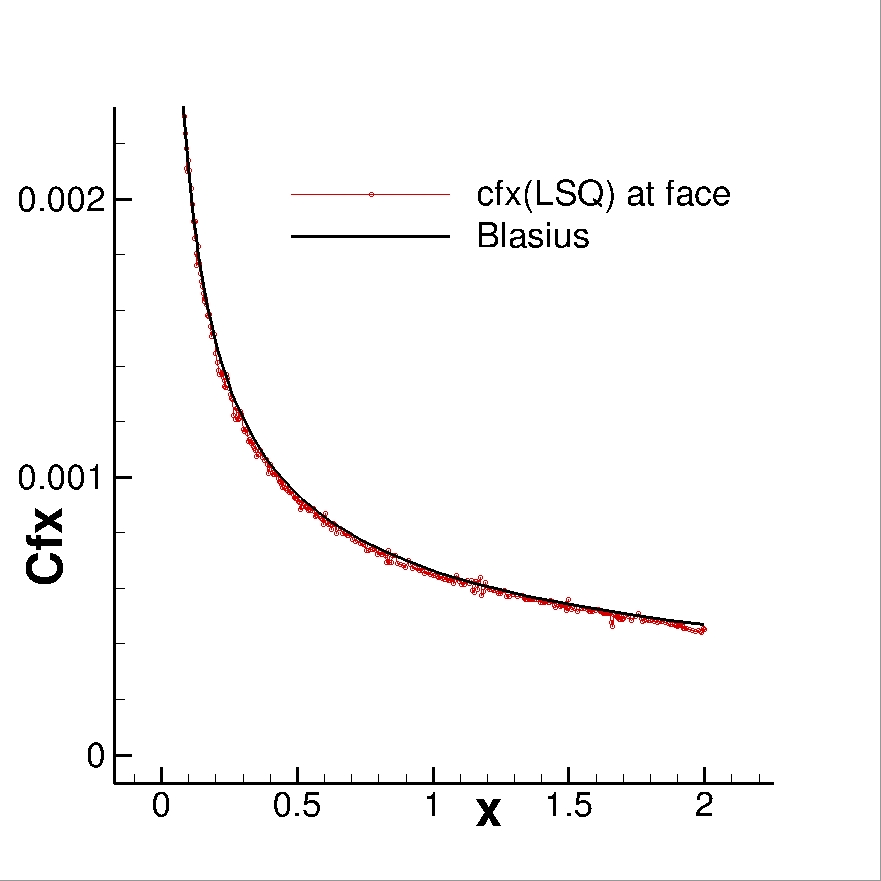}
          \caption{Face-averaged gradients.}
          \label{fig:twod_laminar_fp_re1000000_grid_LSQ_face}
      \end{subfigure}
      \begin{subfigure}[t]{0.32\textwidth}
        \includegraphics[width=\textwidth,trim=4 4 4 4,clip]{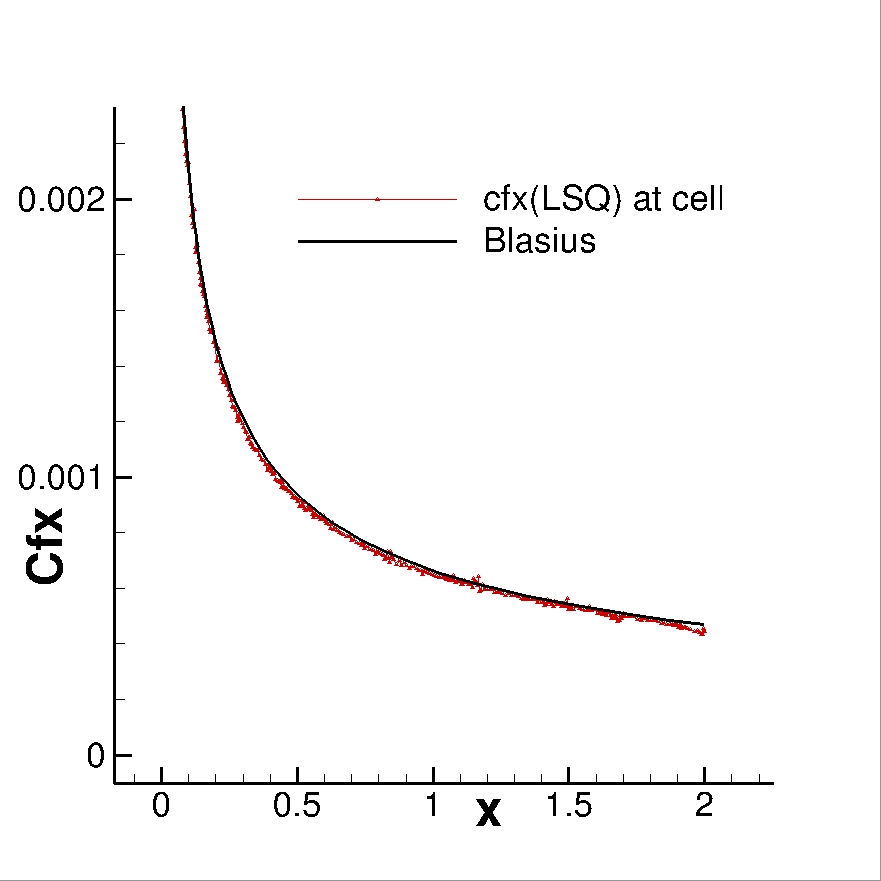}
          \caption{Cell-averaged gradients.}
          \label{fig:twod_laminar_fp_re1000000_grid_LSQ_cell}
      \end{subfigure}
      \caption{Results for a laminar flow over a flat plate at $M_\infty = 0.15$ and $Re_\infty = 10^6$ on an adaptive triangular grid. The flat plate is located at $y=0$ in $x \in [0, 2]$. 
   The skin friction coefficient $C_{fx}$ is plotted at wall face centers, except for the figure (g). }
\label{fig:twod_laminar_fp_re1000000} 
\end{figure}

An adaptive grid was generated by the open-source tool {\it refine} [\citen{park_darmofal:AIAA2008-917}] developed by Mike Park at NASA Langley Research Center and available at \url{https://github.com/nasa/refine}. Adaptation was performed with the multiscale metric technique to control interpolation error in L2-norm [\citen{ugawg-aiaa-verification}] 
with Mach Hessian computed at nodes (by a least-squares method) from the numerical solution obtained at cells with the F-ANG finite-volume solver for the Navier-Stokes equations.
The resulting adaptive grid is shown in Figures \ref{fig:twod_laminar_fp_re1000000_grid} and \ref{fig:twod_laminar_fp_re1000000_grid_zoom}. 
As can be seen, it is a fully irregular grid towards the flat plate located at the bottom in $ x \in [0,2]$. Note that we already have a numerical solution on this grid as 
shown by Mach contours in Figure \ref{fig:twod_laminar_fp_re1000000_grid_zoom}. Then, the skin friction is computed at the flat plate as
\begin{eqnarray}
C_{fx}  
 =  ( \tau_{xx}, \tau_{xy} ) \cdot \hat{\bf n} \frac{2}{M_\infty^2}  = 
\mu \left(  
\frac{2}{3} \frac{\partial u}{ \partial x }  - \frac{1}{3}  \frac{\partial v}{ \partial y}   ,   
 \frac{\partial u}{ \partial y }  +  \frac{\partial v}{ \partial x} 
 \right)  \cdot \hat{\bf n}  \frac{2}{M_\infty^2}
 = 
 \mu \left(   \frac{\partial u}{ \partial y }  +  \frac{\partial v}{ \partial x}  \right) \frac{2}{M_\infty^2},
\end{eqnarray}
where $u$ and $v$ are the $x$- and $y$-components of the velocity, the factor $1/M_\infty^2$ comes from the nondimensionalization (i.e., the velocity
is nondimensionaliezd by the free stream speed of sound), and $\hat{\bf n}$ is either $\hat{\bf n}_i$ or $\hat{\bf n}_j$, depending on where the gradients are computed.
The term $\partial v /  \partial x$ is ignored in the finite-difference methods, but retained in the LSQ gradient methods.

Results are shown in Figures \ref{fig:twod_laminar_fp_re1000000_grid_FD_I}-\ref{fig:twod_laminar_fp_re1000000_grid_LSQ_cell}, where 
numerical skin-friction distributions are compared with that of the Blasius solution. All numerical skin-friction distributions are close to the 
theoretical distribution on average, but significant differences are observed in the noise level. First, the FD-I and FD-II methods 
produce noisy skin-friction distributions as shown in Figures \ref{fig:twod_laminar_fp_re1000000_grid_FD_I} and \ref{fig:twod_laminar_fp_re1000000_grid_FD_II}.
On the other hand, as expected, the proposed methods produce much smoother distributions as shown in Figures \ref{fig:twod_laminar_fp_re1000000_grid_FD_III_eta}
and \ref{fig:twod_laminar_fp_re1000000_grid_FD_III_eta_x}. Note that the choice of $\overline{h}$ has no major impact on the skin-friction distribution: 
Figure \ref{fig:twod_laminar_fp_re1000000_grid_FD_III_eta} for $\eta=0.5$, and Figure \ref{fig:twod_laminar_fp_re1000000_grid_FD_III_eta_x} for
a global constant corresponding to $\eta=0.5$ and $x=1$. These results indicate that it is the uniformness of the step-length, whether in $\eta$ or $y$, that affects
 the smoothness of the skin friction distribution. This is encouraging since the method can be extended to more general problems by developing a method for estimating
 a single value of the step-length $\overline{h}$ in a given physical problem.

For the LSQ gradient methods, the direct use of the nodal gradients leads to a noisy skin-friction distribution as shown in Figure \ref{fig:twod_laminar_fp_re1000000_grid_LSQ}.
A slightly smoother distribution is obtained with the face-averaged nodal gradients as shown in Figure \ref{fig:twod_laminar_fp_re1000000_grid_LSQ_face}.
The smoothest distribution is obtained with the cell-averaged nodal gradients as shown in Figure \ref{fig:twod_laminar_fp_re1000000_grid_LSQ_cell}. 
The results indicate that averaging with more nodal gradients leads to a smoother distribution, but not as smooth as those obtained with the finite-difference III. 
Note that the gradients are linearly exact no mater how they are averaged, and the averaging operator acts as a kind of filtering to remove
high-frequency contents from the underlying gradient approximation.

\section{Concluding Remarks}

We have proposed improved wall-normal derivative formulae for adaptive anisotropic simplex-element grids, based on the simple
idea of unifying the finite-difference distance over a wall. The improved formulae have been demonstrated for a laminar flow over a flat plate 
to reduce the amount of numerical noise in the skin friction distribution. In the future work, these formulae will have to be extended 
to three dimensions and to turbulent flows.


\addcontentsline{toc}{section}{Acknowledgments}
\section*{Acknowledgments}

This work was supported by the Hypersonic Technology Project, through the Hypersonic 
Airbreathing Propulsion Branch of the NASA Langley Research Center, under Contract No. 80LARC17C0004. 
 The author would like to thank Mike Park (NASA Langley Geometry Laboratory) for his assistance with the use of refine. 

\addcontentsline{toc}{section}{References}
\bibliography{../../bibtex_nishikawa_database}

\begin{thebibliography}{10}
\newcommand{\enquote}[1]{``#1''}

\bibitem{WhiteNishikawaBaurle_scitech2020}
White, J., Nishikawa, H., and Baurle, R., \enquote{A 3-{D} Nodal-Averaged
  Gradient Approach for Unstructured-grid Cell-centered Finite-volume Methods
  for Application to Turbulent Hypersonic Flow,} {\em SciTech 2020 Forum\/},
  {AIAA} Paper 2020-0652, Orlando, FL, 2020.

\bibitem{Kleb_etal_aiaa2019-2948}
Kleb, W.~L., Park, M.~A., Wood, W.~A., Bibb, K.~L., Thompson, K.~B., and Gomez,
  R.~J., \enquote{Sketch-to-Solution: An Exploration of Viscous {CFD} with
  Automatic Grids,} {\em 24th {AIAA} Computational Fluid Dynamics
  Conference\/}, {AIAA} Paper 2019-2948, Dallas, TX, 2019.

\bibitem{Nishikawa_scitech2020}
Nishikawa, H., \enquote{A Face-Area-Weighted Centroid Formula for Reducing Grid
  Skewness and Improving Convergence of Edge-Based Solver on Highly-Skewed
  Simplex Grids,} {\em SciTech 2020 Forum\/}, {AIAA} Paper 2020-1786, Orlando,
  FL, 2020.

\bibitem{nishikawa_centroid:JCP2020}
Nishikawa, H., \enquote{A Face-Area-Weighted Centroid Formula for Finite-Volume
  Method That Improves Skewness and Convergence on Triangular Grids,} {\em J.
  Comput. Phys.\/}, Vol.~401, 2020, pp.~109001.

\bibitem{uns_grid_adaptation_aiaa2018-1103}
Park, M.~A., Barral, N., Ibanez, D., Kamenetskiy, D.~S., Krakos, J.~A., Michal,
  T.~R., and Loseille, A., \enquote{Unstructured Grid Adaptation and Solver
  Technology for Turbulent Flows,} {\em 56th {AIAA} Aerospace Sciences
  Meeting\/}, {AIAA} Paper 2018-1103, Kissimmee, Florida, 2018.

\bibitem{alauzet-loseille-decade-aniso-adapt-cfd}
Alauzet, F. and Loseille, A., \enquote{A Decade of Progress on Anisotropic Mesh
  Adaptation for Computational Fluid Dynamics,} {\em Computer-Aided Design\/},
  March 2015, pp.~13--39.

\bibitem{liu_nishikawa_aiaa2016-3969}
Liu, Y. and Nishikawa, H., \enquote{Third-Order Inviscid and Second-Order
  Hyperbolic {N}avier-{S}tokes Solvers for Three-Dimensional Inviscid and
  Viscous Flows,} {\em 46th {AIAA} Fluid Dynamics Conference\/}, {AIAA} Paper
  2016-3969, Washington, D.C., 2016.

\bibitem{nishikawa_hyperbolic_poisson:jcp2020}
Nishikawa, H., \enquote{A hyperbolic {P}oisson solver for tetrahedral grids,}
  {\em J. Comput. Phys.\/}, Vol.~409, May 2020, pp.~109358.

\bibitem{NovelGradStencil:CF2018}
Xiong, M., Deng, X., Gao, X., Dong, Y., Xu, C., and Wang, Z., \enquote{A Novel
  Stencil Selection Method for the Gradient Reconstruction on Unstructured Grid
  Based on {O}pen{FOAM},} {\em Comput. Fluids\/}, Vol.~172, 2018, pp.~426--442.

\bibitem{SozerChristophCetin:AIAA2014}
Sozer, E., Brehm, C., and Kiris, C.~C., \enquote{Gradient Calculation Methods
  on Arbitrary Polyhedral Unstructured Meshes for Cell-Centered CFD Solvers,}
  {\em Proc. of 52nd AIAA Aerospace Sciences Meeting\/}, {AIAA} Paper
  2014-1440, National Harbor, Maryland, 2014.

\bibitem{ShimaKitamuraHaga_AIAAJ2013}
Shima, E., Kitamura, K., and Haga, T.,
  \enquote{{G}reen-{G}auss/Weighted-Least-Squares Hybrid Gradient
  Reconstruction for Arbitrary Polyhedra Unstructured Grids,} {\em {AIAA}
  J.\/}, Vol.~51, No.~11, 2013, pp.~2740--2747.

\bibitem{WangRenPanLi:JCP2017}
Wang, Q., Ren, Y.-X., Pan, J., and Li, W., \enquote{Compact High Order Finite
  Volume Method on Unstructured Grids {III}: Variational Reconstruction,} {\em
  J. Comput. Phys.\/}, Vol.~337, 2017, pp.~1--26.

\bibitem{nishikawa_igg:JCP2018}
Nishikawa, H., \enquote{From Hyperbolic Diffusion Scheme to Gradient Method:
  Implicit {G}reen-{G}auss Gradients for Unstructured Grids,} {\em J. Comput.
  Phys.\/}, Vol.~372, 2018, pp.~126--160.

\bibitem{NishikawaWhite_FANG:jcp2020}
Nishikawan, H. and White, J.~A., \enquote{An efficient cell-centered
  finite-volume method with face-averaged nodal-gradients for triangular
  grids,} {\em J. Comput. Phys.\/}, Vol.~411, 2020, pp.~109423.

\bibitem{NishikawaWhite_scitech2021-xxxx}
White, J.~A. and Nishikawa, H., \enquote{Computation of the Wall Heat Flux for
  a Hypersonic Flow Using Tetrahedral Grids,} {\em AIAA SciTech\/}, 2021, to be
  published as AIAA paper.

\bibitem{Roe_JCP_1981}
Roe, P.~L., \enquote{Approximate {R}iemann Solvers, Parameter Vectors, and
  Difference Schemes,} {\em J. Comput. Phys.\/}, Vol.~43, 1981, pp.~357--372.

\bibitem{nishikawa:AIAA2010}
Nishikawa, H., \enquote{Beyond Interface Gradient: A General Principle for
  Constructing Diffusion Schemes,} {\em Proc. of 40th {AIAA} Fluid Dynamics
  Conference and Exhibit\/}, {AIAA} Paper 2010-5093, Chicago, 2010.

\bibitem{park_darmofal:AIAA2008-917}
Park, M.~A. and Darmofal, D.~L., \enquote{Parallel Anisotropic Tetrahedral
  Adaptation,} {\em Proc. of 46th {AIAA} Aerospace Sciences Meeting and
  Exhibit\/}, {AIAA} Paper 2008-917, Reno, Nevada, 2008.

\bibitem{ugawg-aiaa-verification}
Park, M.~A., Balan, A., Anderson, W.~K., Galbraith, M.~C., Caplan, P.~C.,
  Carson, H.~A., Michal, T., Krakos, J.~A., Kamenetskiy, D.~S., Loseille, A.,
  Alauzet, F., Frazza, L., and Barral, N., \enquote{Verification of
  Unstructured Grid Adaptation Components,} AIAA Paper 2019--1723, 2019.

\end{thebibliography}
\bibliographystyle{aiaa}


\end{document}